\documentclass{amsart}

\makeindex

\usepackage[all]{xy}
\usepackage[english]{babel}

\usepackage{amsmath}

\usepackage{amsthm}

\usepackage{amscd}

\usepackage{amsfonts}

\usepackage{amssymb}

\usepackage{mathrsfs}

\newtheorem{Lem}{Lemma}[section]

    \newtheorem{Prop}[Lem]{Proposition}

    \newtheorem{Thm}[Lem]{Theorem}

    \newtheorem{Cor}[Lem]{Corollary}

\theoremstyle{definition}

   \newtheorem{Def}[Lem]{Definition}

    \newtheorem{Rem}[Lem]{Remark}

\newcommand{\Ps}{\mathbb{P}^}
\newcommand{\Z}{\mathbb{Z}}

\newcommand{\ra}{\rightarrow}
\newcommand{\lra}{\longrightarrow}
\newcommand{\R}{\mathbb R}
\newcommand{\Q}{\mathbb Q}
\newcommand{\C}{\mathbb C}

\begin{document}
\title{Heegner divisors in the moduli space\\ of genus three curves}
\author{Michela Artebani}
\address{Dipartimento di Matematica, Universit\`a di Milano, via C. Saldini 50, Milano, Italia}
\email{michela.artebani@unimi.it}
\thanks{This work is partially supported by: PRIN 2003: Spazi di moduli e teoria di Lie; GNSAGA}
\date{\today}
\begin{abstract}
S. Kond\=o used periods of $K3$ surfaces to prove that the moduli
space of genus three curves is birational to an arithmetic quotient
of a complex 6-ball.
In this paper we study Heegner divisors in the ball quotient, given by arithmetically defined hyperplane sections of the ball.  We show that the
corresponding loci of genus three curves are given by: hyperelliptic
curves, singular plane quartics and plane quartics admitting certain
rational ``splitting curves".
\end{abstract}
\maketitle
\pagestyle{myheadings}
\markboth{Michela Artebani}{Heegner divisors in the moduli space of genus three curves}
\setcounter{tocdepth}{1}
\tableofcontents
\section*{Introduction}
The degree four cyclic cover of the projective plane branched along
a smooth plane quartic is a $K3$ surface endowed with an
automorphism group $G\cong \Z_4$. This simple geometric construction
relates the moduli space $\mathcal M_3$ of genus three curves to a
moduli space $\mathcal M$ of polarized $K3$ surfaces. A generator
$\sigma$ for $G$ can be chosen such that the period point of the
polarized $K3$ surface $(X,\sigma)$ belongs to the $i$-eigenspace
$W_{\C}$ of $\sigma^*$ in $H^2(X,\C)$. This implies that the period
domain is a six dimensional complex ball $B\subset \mathbb
P(W_{\C})$ and the moduli space $\mathcal M$ is obtained by taking
the quotient of $B$ by the action of an arithmetic group $\Gamma$.
The correspondence between genus three curves and polarized $K3$
surfaces is thus described by a period map
$$\mathcal P:\mathcal M_3\backslash \mathcal M_3^h\lra \mathcal M\cong B/\Gamma,$$
where $\mathcal M_3^h$ denotes the hyperelliptic locus.
In fact, in \cite{K} S. Kond\=o proves that this map is birational and induces an isomorphism:
$$\mathcal M_3\backslash \mathcal M_3^h\cong \mathcal M\backslash(\mathcal D_n\cup \mathcal D_h),$$
where $\mathcal D_n$,$\mathcal D_h$ are two irreducible
divisors, called \emph{mirrors}.



The $\Q$-vector space $H^2(X,\Q)$ and the action of $\sigma$ give
naturally a vector space $W$ over $k=\mathbb Q[i]$ such that
$W_{\C}= W\otimes_{k} \C$. An interesting class of divisors in
$\mathcal M$ is then given by quotients of hyperplane sections of
$B$ defined over $k$. These divisors are again arithmetic quotients
of a (five dimensional) complex ball and are called \emph{Heegner
divisors}. This kind of divisors has been introduced and studied in
a more general setting in \cite{Lo}.

The first examples of Heegner divisors are the two mirrors $\mathcal D_n$, $\mathcal D_h$. In \cite{K} it is proved that the generic points in the mirrors
correspond, via the period map $\mathcal P$, to a plane quartic with
a node and to a smooth hyperelliptic genus three curve respectively.
The aim of this paper is to describe all Heegner divisors in terms
of genus three curves. We prove indeed that any Heegner divisor
which is not a mirror can be interpreted as the locus of plane
quartics admitting a \emph{splitting curve} i.e. an irreducible
plane curve whose inverse image by the degree four cyclic cover is
the union of four distinct curves on the $K3$ surface. More
precisely, we show that the splitting curve can be chosen to be rational. The simplest example of Heegner divisor of this kind is
the divisor of plane quartics with a hyperflex line.\\

The paper is organized as follows.

In the first section we introduce some basic notations for genus
three curves, Del Pezzo surfaces and $K3$ surfaces.

The geometric construction by Kond\=o and its main theorem are
recalled in the second section.

In the third section we define Heegner divisors and show how they
are connected to the problem of embedding a rank two lattice in a
non-unimodular lattice. In fact, we introduce two natural invariants
associated to such embeddings: the type $n\in \Z$, $n>0$ and the
index $m\in\{1,2\}$. We also provide an existence result for
Heegner divisors with fixed type and index.

In the following section Heegner divisors are interpreted in terms
of genus three curves. The main theorem states that the generic
point in a Heegner divisor of type $n>1$ corresponds to a smooth
plane quartic having a splitting rational curve with smooth
preimages in the Del Pezzo surface. Moreover, we give the minimal
degree of such a splitting curve as a function of $n$ and $m$ (it
turns out that it is even if $m=1$ and odd if $m=2$).

Some examples of Heegner divisors are provided in the last section
for splitting curves of degree 1, 2 and 3.

The Appendix recalls some results of lattice theory which will be useful in the proofs.\\

\emph{Acknowledgements.} I would like to thank Professor B. van
Geemen for introducing me to this subject and for several valuable
suggestions. I also wish to thank Dr. A. Laface for many interesting
discussions. Part of this work was done during a visit to the
Queen's University (Kingston, Ontario), I thank the Mathematics
Department for hosting me.
\section{Notation and Preliminaries}\label{not}
In this section we introduce the three geometric objects
of main interest in this paper and fix the relative notations: genus three curves, Del Pezzo surfaces and $K3$ surfaces. Our main references are \cite{M}, \cite{DO} and \cite{BPV}.\\

Let $\mathcal M_3$ be the coarse moduli variety parametrizing isomorphism classes of smooth \emph{genus three curves}. We denote with $\mathcal M_3^h$ the hyperelliptic divisor in $\mathcal M_3$.\\

A \emph{Del Pezzo surface} is a smooth surface with ample anti-canonical bundle. These surfaces can be equivalently defined as the blowing up of $\Ps 2$ in a set of $m\leq 8$ distinct points ``in general position".
We are interested in Del Pezzo surfaces of degree two i.e. the blowing up
of $\Ps 2$ in seven points $p_1,\dots, p_7$. In this case, a natural basis for the Picard lattice of the surface $S$ is given by
$$Pic(S)=\langle e_0, e_1, \dots, e_7\rangle,$$
where $e_0$ is the pull-back of the hyperplane bundle of $\Ps2$ and $e_i$ is the exceptional divisor corresponding to $p_i$, $i=1,\dots,7$.
In this basis, the anti-canonical bundle of $S$ is given by:
$$-k=3e_0-(e_1+\cdots+e_7).$$
The morphism $\phi_{(-k)}:S\lra \Ps 2$ associated to the
anti-canonical linear system is a double cover of the plane branched
along a smooth plane quartic curve $Q(S)$. There are exactly 56
effective divisors in $Pic(S)$ with self-intersection $-1$, given by
the preimages of the 28 bitangents of $Q(S)$. Consider now the
lattice:
$$R=k^{\perp}=\{x\in Pic(S): (x,k)=0\},$$
where $(\,,\,)$ denotes the intersection form on $Pic(S)$. This a
root lattice of type $E_7$ and a root basis is given by:
$$\alpha_0=e_0-(e_1+e_2+e_3),$$
$$\alpha_i=e_i-e_{i+1},\ \ \ i=1,\dots,6.$$
We denote with $W(E_7)$ the associated Weyl group.

Let $i$ be the covering involution of $\phi_{(-k)}$, then the induced
involution $i^*$ on $Pic(S)$
generates the center of the Weyl group $W(E_7)$ (see 4., Ch.VII,
\cite{DO}). It is easy to see that
$$i^*(e_i)=-k-e_i,\ \ \ i=1,\dots,7,$$
$$i^*(e_0)=-3k-e_0.$$
The eigenlattices $H_{\pm}$ of $i^*$ relative to $\pm 1$ are given
by
$$H_+=\langle k\rangle,\ H_-=R.$$

Let $X$ be a \emph{$K3$ surface} i.e. a surface with $q(X)=0$ and trivial canonical bundle. We denote with $\omega_X$ a generator for the vector space $H^{2,0}(X)$ of holomorphic two-forms on $X$, with $Pic(X)$ the Picard lattice and with $T(X)$ the transcendental lattice of $X$.
We recall that the cohomology group $H^2(X,\Z)$ is an even unimodular lattice isometric to $L_{K3}=U^{\oplus 3}\oplus E_8^{\oplus 2}$.
The intersection form on $L_{K3}$ is denoted with $(\,,\,)$ .

\section{The model by Kond\=o}\label{K}
\subsection{Geometric construction}
We  briefly recall the geometric construction introduced by Kond\=o
in \cite{K}. Let $C$ be a smooth quartic curve in $\Ps2$ defined by
a homogeneous polynomial $f\in \mathbb C[x,y,z]$ of degree $4$:
$$C=\{(x,y,z)\in\Ps 2:\ f(x,y,z)=0\}.$$
Consider the 4:1 cyclic cover of $\Ps 2$ branched along the divisor
$C$:
$$\pi:X_C\stackrel{4:1}{\lra}\Ps 2.$$
The surface $X_C$ is a $K3$ surface with a degree four polarization, in coordinates it can be given by:
$$X_C=\{(x,y,z,t)\in \Ps 3:\ t^4=f(x,y,z)\}.$$

Let $\sigma$ be a generator for the covering transformation group of $\pi$. We can assume that:
$$\sigma(x,y,z,t)=(x,y,z,it).$$
Let $\tau=\sigma^2$, then the morphism $\pi$ factors naturally through
the double cover of $\Ps 2$ branched along $C$:
$$\pi_1:S_C\stackrel{2:1}{\lra} \Ps 2,$$
where $S_C=X_C/<\tau>$ is a Del Pezzo surface of degree 2 and
$\pi_1=\phi_{(-k)}$ is the morphism associated to the anti-canonical
linear system of $S_C$.
The geometry of the
above construction is then described by the following commutative diagram:
$$\begin{array}{cccccc}
\ \ X_C&\stackrel{\pi_2}{\longrightarrow}&\ \ S_C\\
\ \ \ \ \stackrel{\ \ \ \ \ \ \pi}{\ \searrow} &\ \ \ \ \ &\stackrel{\pi_1\ \ \ }{\swarrow}\\
 & \Ps2&
\end{array}$$

The $K3$ surface $X_C$ is $L_+$-polarized, where $L_+$ denotes the pull-back of $Pic(S_C)$ and the isometry $\sigma^*$ induced by $\sigma$ on the cohomology
lattice acts as a 4-th root of unity on the vector space of holomorphic
two-forms. In other words, the period point of $X_C$ belongs to $\mathbb P(W_{\C})$, where $W_{\C}$ is an eigenspace of $\sigma^*$ in ${L}^{\perp}_+\otimes \C$. In fact, if $L_-=L^{\perp}_+$, then $W_{\C}\cong L_-\otimes_{\Z}\R$ and the action of $\sigma^*$ gives to $W=L_-\otimes_{\Z} \Q$ the structure of a vector space
over $k=\Q[i]$.

The vector space $W_{\C}$ is equipped with the hermitian form
$$\varphi(z,w)=(z,\bar w)$$
of signature $(1,6)$.
The period domain $\mathcal M$ for polarized $K3$ surfaces $(X_C,\sigma^*)$ is
then an arithmetic quotient of the six dimensional complex ball $B$ in
$\mathbb P(W_{\C})$ defined by $\varphi(z,z)<0$:
$$\mathcal M=B/\Gamma,$$
$$\Gamma=\{\gamma\in O(L_-):\gamma\circ\sigma^*=\sigma^*\circ\gamma\}\cong U(\varphi)\cap M(7,\Z[i]).$$

The above construction defines a period map:
$$\mathcal P: \mathcal M_3\backslash\mathcal M_3^h\lra \mathcal M,$$
where $\mathcal M_3^h$ denotes the hyperelliptic locus. In fact, the
image of the period map lies in the complement of two divisors
$\mathcal D_n, \mathcal D_h$, called  \emph{mirrors}, corresponding to non-ample polarized
$K3$ surfaces.
The main result in \cite{K} is:
\begin{Thm}[S. Kond\=o, Theorem 2.5, \cite{K}]
The period map gives an isomorphism:
$$\mathcal P:\mathcal M_3\backslash\mathcal M_3^h\lra \mathcal M\backslash(\mathcal D_n\cup\mathcal D_h).$$
\end{Thm}

\subsection{Isometries}
We now describe the actions of covering transformations on
the cohomology lattice. Let $\tau^*$ be the involution on $H^2(S_C,\mathbb Z)$ induced by $\tau$,
its eigenlattices are given by:
$$L_{\pm}=\{x\in H^2(X_C,\mathbb Z):\tau^*(x)=\pm x\}.$$
In fact, the invariant lattice $L_+$ is the Picard lattice
of the generic $K3$ surface $X_C$ and $L_-=L^{\perp}_+$.
The isomorphism classes of the eigenlattices can be easily computed:
\begin{Lem}\label{iso}
$$L_+\simeq \langle 2\rangle\oplus {A}^{\oplus 7}_1,\ \ L_-\simeq \langle 2\rangle^{\oplus 2}\oplus D_4^{\oplus 3}.$$
\end{Lem}
\proof
From Theorem 4.2.2, \cite{N2} and Theorem 3.6.2, \cite{N1} it follows that $L_+$ is a 2-elementary lattice of signature $(1,7)$ with $r(L_+)=\ell (L_+)=8$ and $\delta(L_+)=1$. Hence obviously: $r(L_-)=14$, $\ell(L_-)=8$ and $\delta(L_-)=1$.
By Theorem \ref{2el} the isomorphism classes of the two lattices are determined uniquely by the set of invariants $(s_+,s_-,\ell,\delta)$.
Hence it is enough to check that the lattices in the right hand sides have the same set of invariants, this follows from the remarks in the Appendix.
For the proof of the first isomorphism see also \S 2, \cite{K}.\qed\\

From Lemma \ref{iso} it follows that $L_+=\pi^*_2 Pic(S_C)$.
For $x\in Pic(S_C)$ we define:
$$\tilde x={\pi}^*_2(x)\in L_+.$$
In particular, fixed a blowing-up morphism $b:S_C\lra \Ps 2$ with the
notation in section \ref{not}, a basis for $L_+$ is given by:
$$\tilde e_0, \tilde e_1,\dots, \tilde e_7.$$
The action of $\sigma^*$ on $L_+$ is clearly the pull-back along
$\pi_2$ of the action of the involution $i^*$ on $Pic(S)$ defined in
the previous section. In particular, let $\tilde R={\pi}^*_2(R)$
and $\widetilde{H}_{\pm}$ be the eigenlattices of $\sigma^*$ on $L_+$.
From the remarks in section 1 immediately follows:
\begin{Lem}
$$\widetilde{H}_+=\langle \tilde k \rangle\cong \langle 4\rangle,\
\widetilde{H}_-=\tilde R\cong E_7(2).$$
\end{Lem}

Moreover, the action of the Weyl group $W(E_7)$ lifts to $L_+$ giving a subgroup:
$$W(E_7)\subset O(L_+)$$
of isometries commuting with $\sigma^*$.\\

A natural basis for $L_-$ is given by:
$$t_1,t_2,f^1_1,\dots,f^1_4,f^2_1,\dots,f^2_4,f^3_1,\dots,f^3_4,$$
where $t_1,t_2$ is a basis for $\langle 2\rangle^{\oplus 2}$ and $f^i_1,\dots,f^i_4$ is a basis for the $i$-th copy of $D_4$ as in the Appendix. Then we have
\begin{Lem}
There exists an isomorphism $L_-\cong \langle 2 \rangle^{\oplus 2}\oplus D^{\oplus 3}_4$ such that the action of the isometry $\sigma^*$ on $L_-$ preserves $\langle 2\rangle^{\oplus 2}$ and each copy of $D_4$. In fact it is given by the matrix $J_1\oplus J^{\oplus 3}_2$ (see the Appendix)
with respect to the previously defined basis of $L_-$.
\end{Lem}
\proof See \cite{ext}.
\qed\\

The discriminant groups of the eigenlattices of $\tau^*$ are given
by (Proposition \ref{es2}):
$$A_{L_+}\cong A_{L_-}\cong \Z_2^{\oplus 8}.$$
In fact, we consider the natural bases:
$$\begin{array}{lcr}
A_{L_+}=\langle\tilde e_0/2,\tilde e_1/2,\dots, \tilde e_7/2\rangle,\\
A_{L_-}=\langle t_1/2,t_2/2,\alpha^i_1,\alpha^i_2,\ i=1,2,3\rangle,
\end{array}$$
where $\alpha^i_1,\alpha^i_2$ generate the discriminant group of the $i$-th copy of
$D_4$ in $L_-$ (their expression in terms of $f^i_1,\dots,f^i_4$ is given in the Appendix). By Proposition \ref{emb}, the primitive embedding of $L_+$ in $L_{K3}$
is determined by an isomorphism
$$\gamma:A_{L_+}\lra A_{L_-}\ \mbox{ with }\ q_{L_-}\circ\gamma=q_{L_+}.$$
\begin{Lem}\label{gamma}
The isomorphism $\gamma$ is given by the matrix
$$M(\gamma)=\left(\begin{array}{cccccccc}
0&1&1&1&1&1&1&1\\
1&0&0&0&0&0&0&0\\
1&0&1&0&0&0&0&0\\
1&0&0&1&0&0&0&0\\
1&0&0&0&1&0&0&0\\
1&0&0&0&0&1&0&0\\
1&0&0&0&0&0&1&0\\
1&0&0&0&0&0&0&1
\end{array}\right)$$
up to the action of $O(A_{L_+})$ (with respect to the previously defined bases).
\end{Lem}
\proof The primitive embedding of $L_+$ in $L_{K3}$ is unique by
Theorem 1.14.4, \cite{N2}. Hence, by Proposition \ref{emb}, the
isomorphism $\gamma$ is unique up to the action of the image of the
natural homomorphism
$$\psi:O(L_+)\lra O(A_{L_+}).$$
In fact, by Theorem 1.14.2 in \cite{N2}, the homomorphism
$\psi$ is surjective. Hence, it is enough to check that $M(\gamma)$ defines an isomorphism preserving the quadratic forms.\qed\\

 The following result describes all vectors in $A_{L_+}$ with
assigned self-intersection with respect to the quadratic form
$q_{L_+}$ and the action of $\sigma^*$ on them:
\begin{Lem}[Lemma 2.1, \cite{K}]\label{vec}
The vectors in $A_{L_+}$ can be divided in four classes accordingly
to their self-intersection:
$$\begin{array}{ll}
C_1=\{x\in A_{L_+}\mid q_{L_+}(x)=-1/2\},\\
C_2=\{x\in A_{L_+}\mid q_{L_+}(x)=1/2\},\\
C_3=\{x\in A_{L_+}\mid q_{L_+}(x)=1\},\\
C_4=\{x\in A_{L_+}\mid q_{L_+}(x)=0\}.
\end{array}$$
\begin{itemize}
\item
An element in $C_1$ is represented uniquely by one of the following
$56$ vectors: $$ \frac{\tilde e_i}{2},\ \frac{\tilde k-\tilde
e_i}{2},\ \ \ 1\leq i\leq 7,$$
$$\frac{\tilde{e_0}-\tilde{e_i}-\tilde{e_j}}{2},\ \frac{\tilde
k-\tilde{e_0}+\tilde e_i+\tilde e_j}{2},\ \ \ 1\leq i<j\leq 7. $$
\item
An element in $C_2$ is represented uniquely by one of the following
$72$ vectors:
$$\frac{\tilde e_0}{2},\ \frac{\tilde k-\tilde{e_0}}{2},$$
$$ \frac{\tilde e_i+\tilde e_j+\tilde e_k}{2}, \ \frac{\tilde k-\tilde e_i-\tilde e_j-\tilde e_k}{2},\ \ \ 1\leq i<j<k\leq 7.$$
\item
An element  in $C_3$ is represented uniquely  by one of the
following $64$ vectors:
$$\frac{\tilde k}{2},\ \frac{\tilde e_i-\tilde e_j}{2},\ \ \ 1\leq i<j\leq 7,$$
$$\frac{\tilde e_0-\tilde e_i-\tilde e_j-\tilde e_k}{2},\ \ \ 1\leq i<j<k\leq 7,$$
$$\frac{2\tilde e_0-\sum_i \tilde e_i+\tilde e_j}{2},\ \ \ 1\leq j\leq 7.$$
\item
An element  in $C_4$ is represented uniquely  by one of the
following $64$ vectors:
$$0,\ \frac{\tilde k+\tilde e_i-\tilde e_j}{2},\ \ \ 1\leq i<j\leq 7,$$
$$\frac{\tilde k+\tilde e_0-\tilde e_i-\tilde e_j-\tilde e_k}{2},\ \ \ 1\leq i<j<k\leq 7,$$
$$\frac{\tilde k+2\tilde e_0-\sum_i \tilde e_i+\tilde e_j}{2},\ \ \ 1\leq j\leq 7.$$
\end{itemize}
The isometry $\sigma^*$ acts as:\\
$$\sigma^*(x)=\left\{\begin{array}{cc}
\tilde k/2-x& \mbox{ if }x\in C_1\cup C_2\\
x & \mbox{ if }x\in C_3\cup C_4.
\end{array}\right.$$
\end{Lem}

Let
$$N_+=\{x\in L^*_+:\ q_{L_+}(x)\in \Z/2\Z\}.$$
Notice that:
$$N_+/L_+=\langle \tilde k/2, \tilde \alpha_1/2,\dots, \tilde\alpha_6/2\rangle.$$
We denote with $O(A_{L_+})$ the group of automorphisms of $A_{L_+}$
respecting the discriminant quadratic form $q_{L_+}$. Then we have
\begin{Lem}\label{orb}
The natural homomorphism $W(E_7)\lra O(A_{L_+})$ is surjective. Moreover, the action of $O(A_{L_+})$ on vectors in Lemma \ref{vec} is
transitive on classes $C_1$ and $C_2$ and has two orbits on classes
$C_3$ and $C_4$.
\end{Lem}
\proof
The first assertion follows from Lemma 2.3, \cite{K}.
The elements of $C_1$ are of the form $\tilde x/2$ where $x$ is the class of a $(-1)$-curve on the Del Pezzo surface $S$. The Weyl group $W(E_7)$ acts transitively on $(-1)$-curves (see for example Lemma 4, Ch. V, \cite{DO}), hence the action of $O(A_{L_+})$ is transitive on $C_1$.\\
Consider now the elements in the class $C_2$:
$$\tilde e_0/2,\ v_{ijk}=(\tilde e_i+\tilde e_j+\tilde e_k)/2.$$
Let $\alpha_{ijk}=\tilde e_0-v_{ijk}$, then:
$$s_{\alpha_{ijk}}(v_{ijk})=3\tilde e_0-2v_{ijk}/2=\tilde e_0/2\ (mod\, L_+).$$
Similarly, two vectors $v_{ijk}$ and $v_{i'j'k'}$ are interchanged by the isometry obtained as composition of the three reflections with respect to $e_i-e_{i'}$, $e_j-e_{j'}$, $e_k-e_{k'}$. This gives the assertion for $C_2$.\\
The element $\tilde k/2$ in the class $C_3$ is preserved by the Weyl group. The other $63$ classes are of the form $\tilde x/2$ where $x$ is a (positive) root of $E_7$. Since $W(E_7)$ acts transitively on roots (see for example \cite{DO}), the action is transitive on $C_3$.\\
The elements in $C_4$ are obtained from those in $C_3$ by adding
$\tilde k/2$, hence the assertion follows easily. \qed

\section{Heegner divisors and embeddings}
\begin{Def}
Let $r\in L_-$ be a primitive vector and let $H_r$ be the
hyperplane section of $B$ orthogonal to $r$:
$$H_r=\{z\in B:\ (z,r)=0\}.$$
The \emph{Heegner divisor} $\mathcal D_r$ associated to $r$ is the image of $H_r$ in $\mathcal M$.
A Heegner divisor $\mathcal D_r$ is of \emph{type} $n$ if $r^2=-2n$.
\end{Def}
\begin{Rem}\ \\
i) For $r\in L_-$ let $z(r)=r-i\sigma^*(r)$, then $z(r)\in W$ and
$$H_r=\{z\in W: \varphi(z,z(r))=0\},$$
where $\varphi$ is the hermitian form on $W$ defined in section 2.
This implies that a Heegner divisor can be equivalently defined as the image in $\mathcal M$ of a hyperplane section of $B$ defined over $k=\mathbb Q[i]$.
In the language of \cite{Lo}, Heegner divisors are quotients of an arithmetically defined arrangement of hyperplanes in the ball.\\
ii) It can be easily seen that $H_r=H_{\sigma^*(r)}$ and $H_r=\emptyset$ if $n\leq 0$.
\end{Rem}
With the notation in the Appendix:
$$\Lambda_r=\langle r,\sigma^*(r)\rangle\cong A_1(n)^{\oplus 2}.$$
Hence, giving a Heegner divisor of type $n\geq 1$ in $\mathcal M$ is equivalent to assign a primitive embedding of the lattice $A_1(n)^{\oplus 2}$ in $L_-$ up to the action of the group $\Gamma$.
Given such an embedding, let ${\Lambda}^{\perp}_r$ be the orthogonal complement of $\Lambda_r$ in
$L_-$ and $P$ be the orthogonal
complement of ${\Lambda}^{\perp}_r$ in $L_{K3}$.
We have an embedding
$$L_+\oplus \Lambda_r \subset P$$
with finite quotient group
$$M_r=P/(L_+\oplus \Lambda_r).$$
We start describing the structure of the group $M_r$.
\begin{Prop}\label{M}Let $r\in L_-$ and $M_r$ as above then:\\
a) if $n$ is even then $M_r\cong \Z_2$ and $(r,L_-)=\Z$;\\
b) if $n$ is odd then:
$$M_r\cong \left\{\begin{array}{lll}
\Z_2\oplus \Z_2 & \mbox{ if }\ &(r,L_-)=2\Z,\\
\Z_2 & \mbox{ if }\ &(r,L_-)=\Z.
\end{array}\right.$$
\end{Prop}
\proof
The group $M_r$ has a natural embedding in $A_{L_+}\oplus A_{\Lambda_r}$. Moreover,
since the embedding of $\Lambda_r$ in $P$ is primitive, the projections on the two factors are isomorphisms (see Appendix):
$$M_r\cong M_1\subset A_{L_+},\ M_r\cong M_2\subset A_{\Lambda_r}$$
such that
$${q_{L_+}}_{\mid M_1}={-q_{\Lambda_r}}_{\mid M_2}.$$
Since
$$A_{L_+}\cong \Z_2^{\oplus 8},\ A_{\Lambda_r}\cong \Z_{2n}^{\oplus
2}$$ it follows that:
$$M_r\cong \Z_2^{\oplus m_r},\ \ 0\leq m_r\leq 2.$$

We first prove that the case $m_r=0$ does not appear i.e. $L_+ \oplus \Lambda_r$ can not be primitively embedded in $L_{K3}$.
This means that for any embedding $L_+\oplus \Lambda_r\hookrightarrow L_{K3}$ there is a primitive element $(x,y)\in  L_+ \oplus  \Lambda_r$ such that $(x,y)/2 \in L_{K3}$.
For this, it suffices to prove that, for some $y \in \Lambda_r$, $y/2$ is not trivial in $A_{L_-}$, since in this case (see \S\,2.2):
$$\gamma^{-1} (y/2)+y/2\in L_{K3}.$$
Let $r=(t,d_1,d_2,d_3)\in L_-$ be a primitive vector with $t\in \langle 2 \rangle^{\oplus 2}$ and $d_i\in D_4$, $i=1,2,3$ according to the direct sum decomposition of $L_-$ given in Lemma \ref{iso}.
From the remarks in the Appendix (in particular Remark \ref{d4}) it follows that we have these cases:\\
i) $d_1, d_2, d_3$ are divisible by two, then $t$ is not divisible by two, hence $r/2=t/2$ is not trivial in $A_{L_-}$;\\
ii) at least one among $d_1, d_2, d_3$ is not divisible by two and not all its coefficients are odd, then $r+\sigma^*(r)/2$ is non-trivial in $A_{L_-}$;\\
iii) at least one among $d_1, d_2, d_3$ is not divisible by two and the non-two divisible $d_i's$ all have odd coefficients, then $r/2$ is in $L_-^*$ and it is non-trivial in $A_{L_-}$.\\
Hence $m_r\in\{1,2\}$.

We now characterize the remaining cases. According to the previous discussion, there are two possibilities:\\
i) $r/2\in A_{L_-}$ i.e. $(r,L_-)=2\Z$, then
$$M_r=\langle r/2,\sigma^*(r)/2\rangle.$$
Since $\Lambda_r$ is primitive in $L_-$, we have that $\sigma^*(r)/2\not=r/2$ in $A_{L_-}$. In particular $n$ is not even, since otherwise $q(r/2)=-n/2\in \Z$ hence $\sigma^*(r)/2=r/2$ in $A_{L_-}$ by Lemma \ref{vec}). Therefore $m_r=2$.\\
ii) $r/2\not\in A_{L_-}$ i.e. $(r,L_-)=\Z$,
then $m_r=1$ and
$$M_r=\langle(r+\sigma^*(r))/2\rangle.$$\qed\\

We call $m_r$ the \emph{index} of the Heegner divisor $\mathcal D_r$, $r\in L_-$.
\begin{Prop}[Existence]\ \\
i) There exist Heegner divisors of type $n$ and index $1$ for every $n\geq 1$;\\
ii) there exist Heegner divisors of type $n$ and index $2$ for every odd $n\geq 1.$\\
(Note that the case of even type and index $2$ is excluded by Proposition \ref{M}).
\end{Prop}
\proof
Consider the following vectors in $L_-$
$$
\begin{array}{lcr}
r_1(1)=f^1_1+f^2_1,\ \ r_1(k)=(k-1)t_1+kf^1_1+f^2_1\ \ (k>1),\\
r_2(k)=kt_1+(k+1)f^1_1,\\
r_3(k)=(k-1)t_1+kf^1_1+f^2_1+f^3_1,\\
r_4(k)=t_1+2kt_2+(k+1)(f^1_1+f^1_2)+k(f^2_1+f^2_2),
\end{array}$$
with respect to the basis of $L_-$ defined in section 2.
It is easy to check that the rank two lattice $\Lambda_i(k)=\langle r_i(k),\sigma^*(r_i(k))\rangle$ is primitive in $L_-$, since the rank two minors of the matrix with rows $r_i(k), \sigma^*(r_i(k))$ have no common factor.
Moreover, notice that
$$r_i(k)^2=\left\{\begin{array}{lcr}
-2(2k)&\mbox{if}& i=1,\ \,\,\,\\
-2(2k+1)& \mbox{if}& i=2,3,\\
-2(4k+1)& \mbox{if}& i=4,\ \,\,\,
\end{array}\right.
$$
in particular the Heegner divisor $\mathcal D_{r_i(k)}$ has even type $n=2k$ for $i=1$ and odd type otherwise.
By Proposition \ref{M} the Heegner divisor $\mathcal D_{r_i(k)}$ has index 1
if $i=1,3$ since $(r_i(k),L_-)=\Z$, hence i) is proved.

For $i=2$ it can be easily checked that the index is 2 iff $k$ is odd i.e. the type $n=2k+1\equiv 3\,(mod\,4).$
For $i=4$ the type is $n=4k+1\equiv 1\,(mod\,4)$ and the index is 2 by Proposition \ref{M}, hence ii) is proved.




\qed\\

The following remark will be useful in the next section
\begin{Lem}\label{c3}
Let $\mathcal D_r$ be a Heegner divisor of even type, then $\tilde k/2\not\in M_r$.
\end{Lem}
\proof
Since $\mathcal D_r$ has index 1 by Proposition \ref{M} we have $M_r=\langle(r+\sigma^*(r))/2\rangle.$
From Lemma \ref{gamma} it follows that the isomorphism $\gamma:A_{L_+}\ra A_{L_-}$ sends $\tilde k/2$ to
$$\gamma(\tilde k/2)=(1/2,1/2,0,\dots,0)$$
(in the usual coordinates for $L_-$). Consider a vector:
$$r=(h_1,h_2,a_1,\dots,a_4,b_1,\dots,b_4,c_1,\dots,c_4),$$
such that $\Lambda_r$ is primitive in $L_-$ and assume that $(r+\sigma^*(r))/2=\gamma(\tilde k/2)$ in $A_{L_-}$.
This means that:\\
i) $h_1$ and $h_2$ are odd,\\
ii) $a_3,b_3$ and $c_3$ are even,\\
iii) $a_1+a_2+a_4$, $b_1+b_2+b_4$ and $c_1+c_2+c_4$ are even.\\
In this case it can be easily checked that $(r,L_-)=2\Z$, giving a contradiction by Proposition \ref{M}.\qed
\section{Heegner divisors and genus $3$ curves}
Let $X_r$ be the generic $K3$ surface with period point in a Heegner divisor $\mathcal D_r$, then its Picard number is equal to $10$ and
$$L_+\oplus \Lambda_r\subset Pic(X_r).$$
In this section we describe the loci of genus three curves corresponding to Heegner divisors of given type and index.
We first recall the special role of Heegner divisors of type 1 in Kond\=o's construction.
\subsection{Type $1$}
\begin{Lem}
A period point in $\mathcal M$ belongs to a Heegner divisor of
type 1 iff the corresponding K3 surface is not the four cyclic cover of the plane branched along a smooth plane quartic.
\end{Lem}
\proof See Theorem 2.5,\cite{K}.\qed
\begin{Prop}
The mirrors $\mathcal D_n$, $\mathcal D_h$ are the unique Heegner divisors of type 1:
\begin{itemize}
\item $\mathcal D_n$ has index 1 and its generic point corresponds to a
plane quartic with a node;
\item $\mathcal D_h$ has index 2 and its generic point corresponds to a
smooth hyperelliptic genus three curve.
\end{itemize}
\end{Prop}
\proof By Lemma 3.3, Theorem 4.3 and Theorem 5.4 in \cite{K}.\qed
\begin{Rem}\ \\
i) In \cite{ext} it is proved that the period map $\mathcal P$ can be defined on the blow-up of the GIT moduli space of semistable plane quartics in one point. This period map gives an isomorphism between $\mathcal D_n$ and the locus of stable singular quartics. Besides, it is proved that the Baily-Borel compactification of $\mathcal D_h$ is isomorphic to the GIT moduli space of semistable sets of eight unordered points in $\Ps 1$ (see also \cite{Kh}).\\
ii) The transcendental lattices of the generic $K3$ surface in $\mathcal D_n$ and $\mathcal D_h$ are given respectively by:
$$\begin{array}{lcr}
T_n\cong U^{\oplus 2}\oplus A_1^{\oplus 8},\\
T_h\cong U(2)^{\oplus 2}\oplus D_8.
\end{array}$$
\end{Rem}
\subsection{Type $n>1$}
\begin{Def}
Let $C$ be a plane quartic and $\pi:X_C\lra \Ps2$ be the four cyclic cover of the
plane branched along $C$. An irreducible plane curve $D$ is a \emph{splitting
curve} for $C$ if the inverse image $\pi^{-1}(D)$ is the union of
four distinct curves.
\end{Def}
It is easy to see that a smooth plane curve $D$ is a splitting curve for $C$ if and only if the restricted cover $\pi_{\mid \pi^{-1}(D)}$ is trivial.
If the curve is not smooth it can happen that it is partially normalized by the cover $\pi$ i.e. the four curves in $\pi^{-1}(D)$ are not isomorphic to $D$ and they can intersect each other outside the ramification curve.
We will see examples of this behaviour in the last section.\\

With $X_r$ as before, we call $\pi:X_r\lra \Ps2$ the four cyclic cover of $\Ps 2$, $C_r$ the corresponding branch quartic and $\pi_2:X_r\lra S_r$ the associated double cover of the Del Pezzo surface $S_r$ (see section \ref{K}).
We start with a simple remark
\begin{Lem}\label{ev}
Let $\mathcal D_r$ be a Heegner divisor of type $n>1$, then $C_r$ admits a splitting curve of even degree.
\end{Lem}
\proof
As noticed before, we have $L_+\oplus \Lambda_r\subset Pic(X_r)$. Consider a class of the form $w=\tilde x+r$ where $\tilde x\in L_+$ with $\tilde x^2\geq 2n-2$. Then $w^2\geq -2$, hence $w$ or $-w$ is effective by the Riemann Roch Theorem. Assume that $w$ is effective and write $w=\sum a_iw_i$, where $a_i$ are positive integers and $w_i$ are irreducible curves.
Then we can assume that $w_1=\tilde x_1+r_1$ with $r_1\not=0$.
The $\sigma^*$-orbit of $w_1$ clearly contains four distinct elements: $w_1,\sigma^*(w_1),\tau^*(w_1)=\tilde x_1-r_1, {(\sigma^*)}^3(w_1)=\sigma^*(\tilde x_1)-\sigma^*(w_1)$. Hence the image $\pi(R)$ of a curve $R$ in the class $w_1$ is a splitting curve for $C$.
The degree of $R$ is given by $\deg(R)=(\tilde x_1,-\tilde k)=2(x_1,-k)$.
\qed
\begin{Thm}\label{heg}
Let $\mathcal D_r$ be a Heegner divisor of type $n>1$ (so $r\in L_-$, $r^2=-2n$)
with index $m_r\in\{1,2\}$. Then there exists a  rational splitting
curve $D$ for $C_r$ with the following properties:
\begin{itemize}
\item[1)]
$$\deg(D)=\left\{\begin{array}{ll}
2(n-1)& \mbox{ if }m_r=1,\\
n-2 & \mbox{ if }m_r=2.
\end{array}\right.$$
\item [2)] The inverse image $\pi_1^{-1}(D)$ in $S_r$ is the union of two smooth rational
curves (hence $\pi^{-1}(D)$ is the union of four smooth rational
curves in $X_r$).
\end{itemize}
\end{Thm}
\proof
The idea is to construct the class of a $(-2)$-curve in $Pic(X_r)\backslash L_+\oplus \Lambda_r$ starting from the information given by the structure of the group $M_r$.\\
i) If $m_r=1$ (i.e. $M_r\cong\Z_2$) then there exists $\tilde x\in L_+$ such that:
$$(\tilde x+f)/2\in Pic(X_r)$$
where $f=r+\sigma^*(r)$. Notice that
$(\tilde x'+f)/2\in Pic(X_r)$ for any $\tilde x'\in L^*_+$ with $\tilde x'/2=\tilde x/2$ in $A_{L_+}$.
Moreover
$$q(\tilde x/2)=-q(f/2)=n.$$
Then $\tilde x/2\in A_{L_+}$ belongs to the class $C_3$ in
Lemma \ref{vec} if $n$ is odd and to the class $C_4$ if $n$ is
even.

\noindent If $n$ is odd we can assume (up to the action of the Weyl group, see Lemma \ref{orb} and Lemma \ref{c3}) that:
$$\tilde x/2=(\tilde e_1-\tilde e_2)/2\in A_{L_+}.$$
Let $x'=(n-1)e_0-(n-2)e_1-e_2\in Pic(S_r)$, then $\tilde x'/2=\tilde x/2$ in $A_{L_+}$. Hence by a previous remark
$$y=(\tilde x'+f)/2\in Pic(X_r)$$
and $y^2=(x')^2/2-n=-2.$
Notice that $y+\tau^*(y)=x'$ and the arithmetic genus $p_a(x')$ is zero by the adjunction formula.

The class $y$ can be written as $y=\sum_{i=1}^l a_iy_i$, where $a_i$ are positive integers and $y_i$ are irreducible curves.
Since $p_a(x')=0$, there exists an irreducible component $y_1$ of $y$ which is a smooth rational curve and doesn't belong to $L_+$.
This can be written as
$$y_1=(\tilde x_1+f_1)/2,\ \ \tilde x_1\in L_+,\ f_1\in
\Lambda_r,\ f_1\not=0.$$
Since $y^2=y_1^2=-2$ and $f_1^2\leq f^2$ we have $x_1^2\geq (x')^2=2(n-2)$. Moreover, since $x_1$ is a component of $x'$, its arithmetic genus is zero.
This gives
$$0=p_a(x_1)=\frac{x_1^2+(x_1,k)}{2}+1\geq \frac{2(n-2)+(x_1,k)}{2}+1,$$
hence
$$(y_1,-k)=(x_1,-k)\geq 2(n-1)=(y,-k).$$
This implies that $y_1=y$ i.e. $y$ is a $(-2)$-curve.
Let $R$ be the curve in the class $y$, then the plane curve $D=\pi(R)$ is a rational splitting curve for $C_r$ and:
$$\deg(D)=-(x',k)=2(n-1).$$

If $n$ is even, then we can assume (up to the action of the Weyl group, see Lemma \ref{orb}) that:
$$\tilde x/2=(\tilde k+2\tilde e_0-\sum_i\tilde e_i+\tilde e_1)/2=(\tilde e_0-\tilde e_1)/2\in A_{L_+}.$$
The proof follows as in the odd case with the choice
$x'=(n-1)e_0-e_1-(n-2)e_2$.\\

\noindent ii) If $m_r=2$ (i.e. $M_r\cong \Z_2\oplus \Z_2$) then there exists $\tilde x\in L_+$ such that:
$$(\tilde  x+r)/2,\ \sigma^*(\tilde x+r)/2\in Pic(X_r),$$
where
$$q(\tilde x/2)=-q(r/2)=n/2.$$
Then $\tilde x/2\in A_{L_+}$ belongs to the class $C_1$ in
Lemma \ref{vec} if $n\equiv 3\, (mod\,4)$ and to the class $C_2$ if $n\equiv 1\, (mod\,4)$.

If $n\equiv 3\, (mod\,4)$, then we can assume that:
$$\tilde x/2=\tilde e_1/2.$$
Let $n=4s+3$ and consider the vector $x'=2se_0-(2s-1)e_1.$ Note that
$p_a(x')=0$ by the adjunction formula. Then
$$y=(x'+r)/2\in Pic(X_r)$$
and $y^2=-2$.
It can be proved as in i) that $y$ is a smooth $(-2)$-curve
and if $R$ is a curve representing $y$, then $D=\pi(R)$ is a rational splitting curve with
$$\deg(D)=-(x',-k)=n-2.$$
If $n\equiv 1\, (mod\,4)$, then we can assume that:
$$\tilde x/2=\tilde e_0/2.$$
Let $n=4s+1$ and consider the class $x'=(s-1)(e_0-e_1).$
Defining $y$ as above, similar computations give a rational splitting curve $D$ with:
$$\deg(D)=-(x', k)=n-2.$$\qed
\begin{Rem}
Note that, in the proof of Theorem \ref{heg}, (the classes of) the preimages of the splitting curves in the Del Pezzo surface are given explicitly up to the action of the Weyl group.
\end{Rem}
\begin{Cor}\label{gen}
The Picard lattice of the $K3$ surface $X_r$ is generated by $L_+$ and two smooth rational curves $R, \sigma^*(R)$ with $R\not\in L_+\oplus \Lambda_r$.
\end{Cor}
\proof
Obviously $Pic(X_r)\supset \langle L_+, R, \sigma^*(R)\rangle$.
Let $R=(\tilde x+y)/2$, where $\tilde x\in L_+$ and $y\in \Lambda_r$ are primitive.
Then the thesis follows since the embedding of $L_+$ in $Pic(X_r)$ is primitive and the classes $R,\sigma^*(R)$ are not divisible in $L_{K3}$.\qed\\

We also give a partial converse to the previous result:
\begin{Prop}\label{hegconv}
Let $C$ be a plane quartic admitting a rational splitting curve $D$ such that
the inverse image of $D$ in $S_C$ is the union of two smooth
rational curves with primitive classes in $Pic(S_C)$.
Let $d$ be the minimal degree of a curve $D$ with this
property. Then $X_C$ belongs to a Heegner divisor\\
i) of index $m=2$ and type $n=d+2$ if $d$ is odd;\\
ii) of index $m=1$ and type $n=(d+2)/2$ if $d$ is even.
\end{Prop}
\proof
Let $R$ be one of the inverse images of $D$ in $X_C$.
Notice that:
$$R=(x+y)/2\in Pic(X_C),$$
where
$$x=R+\tau^*(R)\in L_+,\ y=R-\tau^*(R)\in L_-\cap Pic(X_C).$$
An easy computation gives:
$$x^2=2(d-2),\ y^2=-2(d+2).$$
Notice that $y\not=0$ in $L_-$ since otherwise $x=2R$ (which is the class of a preimage of $D$ in $S_C$) would be divisible in $L_+$.
Thus $X_C$ belongs to a Heegner divisor $\mathcal D_r$ with $y\in \Lambda_r$.
Let $\Lambda^{\perp}_r$ be the orthogonal complement of $\Lambda_r$ in $L_-$ and $P_r$ be the orthogonal complement of $\Lambda^{\perp}_r$ in $L_{K3}$. Then we define
$$M_r=P_r/L_+\oplus \Lambda_r.$$
Notice that:
$$\sum_{i=0}^{3} {({\sigma^*})^i(R)}=-d\tilde k.$$
Then we have the following cases:\\
i) If $d$ is odd, then:
$$\sum_{i=0}^3 ({\sigma^*})^i(R)/2=(x+\sigma^*(x))/2=\tilde k/ 2\not= 0 \in A_{L_+}.$$
Hence $x/2\not=\sigma^*(x)/2\in A_{L_+}.$
Therefore
$$M_r\cong \langle x/2, \sigma^*(x)/2\rangle\subset A_{L_+},$$
this means that $\mathcal D_r$ has index $m_r=2$.
By Theorem \ref{heg} (see also the proof) the curve $C$ admits a rational splitting curve $D'$ of degree $d'=-r^2/2-2$ with smooth and primitive preimages in the Del Pezzo surface such that
$\pi^*(D')=\sum_{i=1}^{4} {(\sigma^*)}^i(R'),\ i=1,\dots,4$ with
$R'-\tau^*(R')=r$.
Then:
$$(R'-\tau^*(R'))^2=-2(d'+2)$$
is maximal i.e. $d'$ is minimal.
Then $d=d'$ and $\mathcal D_r$ has type $n=d+2$.\\
ii) If $d$ is even, then:
$$\sum_i {(\sigma^*)}^i(R)/2=(x+\sigma^*(x))/2= 0 \in A_{L_+}.$$
Hence:
$$x/2=\sigma^*(x)/2\in A_{L_+}.$$
The case $x/2= 0$ in $A_{L_+}$ (i.e. $x=2x'\in L_+$) can be excluded since otherwise the inverse images of $D$ in $S_C$ would belong to divisible classes.
Then the Heegner divisor $\mathcal D_r$ has index $m_r=1$:
$$M_r=\langle x/2\rangle\subset A_{L_+}.$$
By applying Theorem \ref{heg} as in i), the minimality condition on the degree gives that
$\mathcal D_r$ has type $n=(d+2)/2$.
\qed
\begin{Cor}\label{odd}
A plane quartic $C$ admits a splitting curve of odd degree if and only if $X_C$ belongs to a Heegner divisor of index $2$.
\end{Cor}
\proof
The result follows from Theorem \ref{heg} and the proof of Proposition \ref{hegconv} if one notices that in the odd degree case the inverse images of the splitting curve in $Pic(S_C)$ always have primitive classes.
\qed
\begin{Rem}\ \\
i) Splitting curves of even degree do not always give rise to Heegner divisors. From the proof of Proposition \ref{hegconv} it follows that this happens iff a inverse image $R$ of the splitting curve in the $K3$ surface satisfies:
$$R-\tau^*(R)=0.$$
For example, let $C$ be the generic plane quartic with equation $f_4(x,y,z)=0$. It easy to prove that a plane quartic $D$ of the form $l(x,y,z)^4=f_4(x,y,z)$, where $l$ is linear,  is a splitting curve for $C$. The inverse image of $D$ in $X_C$ splits in four genus three curves ${(\sigma^*)}^i(R)$, $i=0,\dots,3 $ which are hyperplane sections of $X_C$ in $\Ps3$ (i.e. they all belong to the class $-\tilde k\in L_+$).\\
ii) If the splitting curve has even degree $d\not\equiv 2\,(mod\,4)$ an easy computation shows that its inverse images have primitive classes in the Del Pezzo surface.\\
iii) Theorem \ref{heg} and Proposition \ref{hegconv} describe the loci of plane quartics corresponding to the union of all Heegner divisors of given type and index. However, it is not clear if these loci are irreducible. In other words, it is unknown if an Heegner divisor is determined uniquely by its type and index.
Note that this problem is equivalent to classify the embeddings of the lattice $\Lambda_n$ in $L_-$ up to the action of $\Gamma$.\\
iv) In \cite{Kh} Heegner divisors in $\mathcal D_h$ (in fact, this is a union of our Heegner divisors with fixed type and index intersected with $\mathcal D_h$) are introduced and the existence of an automorphic form vanishing exactly on them is proved.
\end{Rem}
\section{Examples}
In this section we give some examples of Heegner divisors corresponding to plane quartics with splitting curves of degree one, two and three.
We will usually identify Heegner divisors in $\mathcal M$ and their inverse image by $\mathcal P$ in $\mathcal M_3$.
\\

Let $C\subset \Ps2$ be a plane curve of degree $dd'$ with $d=2,4$ and
let $\pi_{d}:Z\ra \Ps2$ be the $d$:1 cover of $\Ps2$ branched along
$C$. We start recalling the following result:
\begin{Prop}[Proposition 1.7, Ch.3, \cite{Ve}]\label{split}
Let $D\subset \Ps2$ be a curve not containing components of $C$. Let
$\tilde D=\pi_{d}^{-1}(D)$. The restriction of the cover:
$${\pi_d}_{|\tilde D}:\tilde D\lra D$$
is trivial if and only if there exists a curve $B\subset \Ps2$ of
degree $d'$ such that:
$$C\cdot D=d B\cdot D.$$
\end{Prop}
By the remark after Definition \ref{split}, this result characterizes smooth splitting curves of $C$.
\subsection{Quartics with a hyperflex}
Let $C$ be a generic smooth plane quartic with a hyperflex line $L$ i.e.
$$L\cdot C=4p,\ \ p\in C.$$
After a projective transformation we can assume that the line $L$ is
given by the equation $x=0$ and that the quartic $C$ is of the
form:
\begin{equation}C:\ xf_3(x,y,z)+z^4=0,\end{equation}
where $f_3$ is a cubic polynomial in $\mathbb C[x,y,z]$.
It is known that the locus of plane quartics with at least one hyperflex line is an irreducible closed subvariety of codimension one in the moduli space of genus three
curves (see Proposition 4.9, Ch.I, \cite{Ve}).
As usual, let $\pi_1:S_C\ra \Ps 2$ and $\pi:X_C\ra \Ps2$ be the double and the 4:1 cover of $\Ps2$ branched along $C$ respectively.
\begin{Prop}
The locus of plane quartics with a hyperflex line is the unique Heegner divisor of type $3$ and index $2$.
\end{Prop}
\proof
Let $C$ be a plane quartic with a hyperflex line $L$. By Proposition \ref{split} both $\pi_1$ and $\pi$ are trivial over $L$. Hence the line $L$ is a splitting curve for $C$, moreover it splits in the union of two smooth rational curves on the Del Pezzo surface $S_C$. By Proposition \ref{hegconv} it follows that $X_C$ belongs to a Heegner divisor of type $n=3$ and index $2$ (the degree of $L$ is obviously minimal). The converse follows from Theorem \ref{heg}.
Since the locus of plane quartics with a hyperflex line is irreducible in $\mathcal M_3$, there is only one Heegner divisor of type $3$ and index $2$.
\qed\\

\noindent We denote this Heegner divisor with $\mathcal D_{flex}$. We now compute the Picard lattice of the generic $K3$ surface $X$
in $\mathcal D_{flex}$.
Let $M,M'$ be the inverse images of $L$ in the Del Pezzo surface $S$.
$${\pi_1}^*(L)=M+M'.$$
Each of them splits in the union of two
smooth rational curves on $X$:
$${\pi_2}^*(M)=M_1+M_2,\ {\pi_2}^*(M')=M_1'+M_2'.$$
Notice that $\sigma^*(M_i)=M'_i$, $i=1,2$, in particular:
$$M_1+M_2,\ M'_1+M'_2\in L_+,$$
$$M_1-M_2,\ M'_1-M'_2\in L_-\cap Pic(X).$$
Let $r=M_1-M_2$, it is easy to check that $r^2=-6.$
Hence the period point of the $K3$ surface $X$ lies in the
hyperplane section of $B$ determined by the vector
$r\in L_-$.
By Corollary \ref{gen} we have:
$$Pic(X)=\langle L_+,M_1,M'_1\rangle.$$
Notice that $M, M'$ are two $(-1)$-curves on $S$ with:
$$M+M'={{\pi_1}^*}(L)=-k.$$
Hence we can assume, up to the action of the Weyl group, that:
$$M=e_7,\ \ M'=3e_0-(e_1+\cdots+e_6+2e_7).$$
Then:
$$M_1+M_2=\tilde e_7,\ \ M'_1+M'_2=3\tilde e_0-(\tilde e_1+\cdots+\tilde e_6+2\tilde e_7).$$
The intersection matrix of the Picard lattice with respect to the basis $\tilde e_0, \tilde e_1,\dots,\tilde e_7,$ $M_1,M'_1$ is then given by:
$$Q(X)=\left(\begin{array}{cccccccc}
2  & O_{1,6}           & 0 & 0   & 3 \\
O_{6,1}  & -2I_{6}& O_{6} & O_{6}   & I_{6,1} \\
0  & 0            &-2 & -1  & 2 \\
0  & 0            &-1 & -2  & 1\\
3  & 1            &2  & 2   & 1\\
\end{array}\right)$$
where $O_{1,6}$ and $O_{6,1}$ are respectively a row and a column
zero vector of lenght $6$, $O_6$ is the zero matrix of order $6$,
$I_{6}$ is the identity matrix of order $6$ and $I_{6,1}$ a column
vector of 1's.

\subsection{Quartics with splitting conics}
We consider the locus of smooth plane quartics $C$
such that there exists an irreducible conic $T$ with intersection
divisor of the form:
$$T\cdot C=4p+4q,$$
where $p\not=q$.
After a projective transformation we can assume that the conic $T$
is given by:
$$T:\ xy-z^2=0$$
and that $p=(0,1,0)$, $ q=(1,0,0).$
Then, up to projectivities, the equation of the quartic $C$ is of
the form:
$$C:\ (xy-z^2)f_2(x,y,z)+z^4=0,$$
where $f_2$ is a quadratic polynomial in $\mathbb C[x,y,z]$.
It is known (see \cite{IH}) that a smooth plane quartic $C$ has $63$ one dimensional families of tangent conics and that each family contains six reducible conics (i.e. the union of two bitangents).
We now prove that
\begin{Lem}\label{ir}
The locus of plane quartics with an irreducible splitting conic is a codimension one closed irreducible suvariety of $\mathcal M_3$.
\end{Lem}
\proof
Let $\Delta\subset \mathbb P^{14}$ be the discriminant locus
corresponding to singular plane quartics. Notice that
polynomials like the one defining $C$ define a $6$ dimensional
irreducible subvariety $A$ of $\mathbb P^{14}\backslash\Delta$. Let
$\phi:\mathbb P^{14}\backslash\Delta\rightarrow \mathcal M_3\backslash\mathcal M^h_3$ be the
natural morphism to the moduli space of non-hyperelliptic genus
three curves. This is a closed and surjective morphism (see
Proposition 4.7, \cite{Ve}). In particular, $V=\phi(A)$ is an
irreducible closed subvariety of $\mathcal M_3\backslash\mathcal M^h_3$. Note that the dimension of the fibers of $\phi$ over $V$ is equal to one. In fact, it equals the dimension of:
$$\{g\in PGL(3,\mathbb C):g(p)=p,\ g(q)=q,\ g(T)=T\}.$$
Thus $V$ has codimension one in $\mathcal M_3\backslash\mathcal M^h_3$.
In fact, it can be proved as in \cite{Ve} that $V$ is also closed in $\mathcal M_3$.
\qed
\begin{Prop}
The locus of plane quartics with an irreducible splitting conic is the unique Heegner divisor of type $2$ and index $1$.
\end{Prop}
\proof By Proposition \ref{split}, the inverse
image of the conic $T$ splits in the union of two smooth rational
curves $R,R'$ on the Del Pezzo surface and in four smooth rational curves $R_1, R_2, R'_1=\sigma^*(R_1), R'_2=\sigma^*(R_2)$ on the $K3$ surface $X_C$.
By \cite{IH}, $R$ is equivalent to the union of two $(-1)$-curves $L_1,M_1$ with $(L_1,M_1)=1$ and mapping to distinct bitangents of $C$.
Notice that $(L_1,\sigma^*(M_1))=0$.
We can assume, up to the action of the Weyl group, that:
$$L_1=e_0-e_1-e_2,\ \sigma^*(M_1)=2e_0-e_1-e_2-e_5-e_6-e_7.$$
Then
$$\begin{array}{lcr}
R=L_1+M_1=2e_0-(e_1+\cdots+e_4),\\
R'=\sigma^*(L_1+M_1)=4e_0-(e_1+\cdots+e_4)-2(e_5\cdots+e_7).
\end{array}$$
In particular, $R$ is a primitive class in $Pic(S_C)$. Then, by Proposition \ref{hegconv}, $X_C$ belongs to a Heegner divisor of type $2$ and index $1$. The converse follows from Theorem \ref{heg} and Lemma \ref{ir} implies that there is a unique Heegner divisor of type $2$ and index $1$.\qed\\

We denote this divisor with $\mathcal D_{con}$. By Corollary \ref{gen} the Picard lattice of the generic $K3$ surface $X$ in $\mathcal D_{con}$ is given by:
$$Pic(X)=<L_+,R_1,R'_1>.$$
The intersection matrix of the lattice $Pic(X)$ with
respect to the basis $\tilde e_0,\tilde e_1,\dots, \tilde e_7,$ $R_1,
R'_1$:
$$Q(X)=\left(
\begin{array}{cccccccccccccc}
 2 &0 &0 &0 &0 &0 &0 &0 &2 & 4\\
0 &-2 &0 &0 & 0 &0 &0 &0 &1 &1\\
0  &0  &-2& 0 &0 &0 &0 &0 & 1 & 1\\
 0  &0  &0& -2 &0 &0 &0 & 0& 1 & 1\\
0  &0  &0& 0 &-2 &0 &0 &0 & 1 & 1 \\
 0  &0& 0 &0 &0 &-2 &0 & 0 & 0 & 2\\
 0  &0& 0 &0 &0 &0 &-2 & 0 & 0 & 2\\
 0  &0& 0 &0 &0 &0 &0 & -2 & 0 & 2\\
  2  &1 &1 &1 &1 & 0 &0 &0 &-2 & 2\\
  4  & 1 &1  &1  &1 & 2 &2 &2 &2 &-2\\
\end{array}
\right)$$
\subsection{Quartics with splitting cubics}
A simple remark is the following:
\begin{Lem}\label{cub}
A plane quartic $C$ has a hyperflex line if and only if there exists
a cubic $D$ such that
the 4:1 cover of $\Ps2$ branched along $C$ is trivial over
$D$.
\end{Lem}
\proof
Notice that a quartic $C$ has a hyperflex line $L$ if and only if there
is a point $p\in C$ with
$$4p\in \mid O_C(1)\mid=\mid K_C\mid.$$
In this case, let $M$ be a line through $p$, $M\not=L$, then $M\cdot C-p=p_1+p_2+p_3.$
Hence
$$4(p_1+p_2+p_3)\in \mid 4K_C-4p\mid=\mid 3K_C\mid$$
i.e. there exists a cubic curve $D$ intersecting $C$ with multiplicity $4$ in each point and such that the points in $C\cdot D$ lie on the line $M$.
Equivalently,  by Proposition \ref{split}, the 4:1 cover of $\Ps2$
branched along $C$ is trivial over $D$.

Conversely, if there exists a cubic $D$ with $D\cdot C=4(p_1+p_2+p_3)$ and a line $M$ with $M\cdot D=p_1+p_2+p_3$ then
$$M\cdot C=p_1+p_2+p_3+p,$$
hence $4p=4M\cdot C-D\cdot C\in  \mid O_C(1)\mid$ i.e. $p$ is a hyperflex point.
\qed
\begin{Rem}
In fact, Lemma \ref{cub} follows immediately if we recall that, up
to a projective transformation, the equation of a quartic $C$ with a
hyperflex line is of the form:
$$xf_3(x,y,z)+z^4=0,$$
where $f_3$ is a cubic polynomial in $\mathbb C[x,y,z]$.
\end{Rem}
In general, a (singular) plane cubic $D$ can be a splitting curve for a quartic $C$ even if the 4:1 cover branched along $C$ is not trivial over $D$ i.e. the irreducible components of $\pi^{-1}(D)$ are not isomorphic to $D$.
Let $\pi:X_C\ra\Ps2$ and $\pi_1:S_C\ra \Ps2$ be the 4:1 and the double cover branched along $C$ respectively.
\begin{Prop}
Let $C$ be a plane quartic with a splitting nodal cubic $D$. Then $X_C$ belongs to a Heegner divisor of index 2 and\\
i) type $3$ if $\pi$ is trivial over $D$ (i.e. $X_C$ lies in $\mathcal D_{flex}$),\\
ii) type $7$ if $\pi_1$ is trivial over $D$,\\
iii) type $5$ otherwise.
\end{Prop}
\proof
Up to projectivities we can assume that $D$ is defined by the equation:
$$y^2=x^3+x^2.$$
A parametrization for $D$ is given by:
$$\phi:\C\lra D,\ \ t\longmapsto(t^2-1,t(t^2-1)),$$
where $\phi(1)=\phi(-1)$ is the node of $D$.

The plane quartic $C$ with equation $f_4(x,y,z)=0$ intersects $D$ in
the points $\phi(t)$, with:
$$f_4(t^2-1,t(t^2-1),1)=0.$$
In particular, if the intersection multiplicity of each point in
$C\cap D$ is a multiple of four, then:
$$f_4(t^2-1,t(t^2-1),1)=q^4,$$
where $q\in \C[t]$ and $q^4\in \C[t^2-1,t(t^2-1)]$.
It can be easily proved that
$$\C[t^2-1,t(t^2-1)]=\{p\in\C[t]:\ p(1)=p(-1)\}.$$
Hence:
$$q(1)=i^aq(-1),\ \ a=0,\dots,3.$$
The 4:1 cover of the plane branched along $C$ is given by:
$$X_C=\{w^4=F(x,y,z)\}\subset\Ps3.$$
Notice that the map $\phi$ lifts to four distinct maps to $X_C$:
$$\phi_k:\C\lra X_C,\ \ t\mapsto(t^2-1:t(t^2-1):1:i^kq(t)),$$
hence the inverse image of $D$ in $X_C$ has four isomorphic
irreducible components:
$$D_k=\phi_k(\C),\ \  a=0,\dots,3.$$
We now study the restriction of the cover $\pi$ to $D$ for all values of $a$.\\
i) If $a=0$, then $q(1)=q(-1)$, hence each component $D_k$ is a nodal curve and $\pi$ is trivial over $D$. By Lemma
\ref{cub} the plane quartic $C$ has a hyperflex line and
$X_C$ belongs to the Heegner divisor $\mathcal D_{flex}$.
\\
ii) If $a=2$, then $q(1)=-q(-1)$, hence the curves $D_k$ are smooth rational curves. Note that $\phi_0(1)=\phi_2(-1)$ and $\phi_0(-1)=\phi_2(1)$, so $D_0$ and $D_2$ meet in two points over the singular point $p=(0:0:1)\in D$. Similarly, $D_1$ and $D_3$ meet in two points over $p$. In this case the double cover $\pi_1$ of $\Ps 2$ branched along $C$ is trivial over $D$, in particular the inverse image of $D$ in $S_C$ is the union of two
singular curves. Note that:
$$r=D_0-\tau^*(D_0)=D_0-D_2\in Pic(X_C)\cap L_-$$
and $r^2=-4-2\cdot 2-2\cdot 3=-14.$
We have a lattice embedding:
$$\Lambda=\langle r,\sigma^*(r)\rangle\subset Pic(X_C)\cap L_-\cong A_1(n)^{\oplus 2}.$$
In particular $-14=-2n(a^2+b^2)$, $a,b\in \Z$. Thus we have $n=7$
and $X_C$ belongs to a Heegner divisor of type $7$ and index $2$.\\
iii) If $a=1$, then $q(1)=iq(-1)$, hence the curves $D_k$ are smooth. In this case we have also $\phi_0(1)=\phi_1(-1)$ and $\phi_0(-1)=\phi_3(1)$, so $D_0$ meets $D_1$ and $D_3$. Similarly, $D_2$ meets $D_3$ and $D_1$.
In particular, the inverse image of $D$ in $S_C$ is the union of
two smooth rational curves. Note that:
$$r=D_0-\tau^*(D_0)=D_0-D_2\in Pic(X_C)\cap L_-$$
and
$$r^2=-4-2\cdot 3=-10.$$
Then:
$$\Lambda=\langle r,\sigma^*(r)\rangle\subset Pic(X_C)\cap L_-.$$
In fact (since $-r^2/2$ is prime and $n\not=1$) we easily get that
this is an equality. Hence $X_C$ belongs to a Heegner divisor of type $5$ and index $2$. This could be proved also by applying Proposition
\ref{hegconv} (note that the degree is minimal since $m=2$ and we
can assume that $C$ has no hyperflex lines).\\
iv) The case $a=3$ is analogous to case iii).
\qed
\begin{Rem}
A similar description could be given more generally for any rational
splitting curve.
\end{Rem}

\appendix
\section{Lattices}
A \emph{lattice} $L$ is a free $\Z$-module of finite rank endowed with an integral symmetric bilinear form
$$(\, , \,):L\times L\lra \Z.$$
Let $r(L)$ be its \emph{rank}, $d(L)$ its  \emph{discriminant} and  $s(L)=(s_+,s_-)$ its \emph{signature}.
We denote with $O(L)$ the group of self-isometries of $L$.

If $L_1, L_2$ are two lattices then $L_1\oplus L_2$ is their
\emph{orthogonal direct sum} and $L^{\oplus n}$ is the orthogonal
direct sum of $n$ copies of $L$. Let $a$ be an integer, then $L(a)$ is the lattice with the same underlying $\Z$-module of $L$ and with bilinear form multiplied by $a$.\\

In the following, we assume that $L$ is an even lattice i.e. $(x,x)$ is even for each $x\in L$.
Let $L^*=Hom(L,\Z)$ be the dual lattice of $L$, then there is a canonical
embedding $L\subset L^*$. The \emph{discriminant group} of $L$ is
the quotient group
$$A_L:=L^*/L.$$
The group $A_L$ is an abelian group of order $\mid d(L)\mid$, we
denote by $\ell(L)$ the minimal number of its generators. The
bilinear form on $L$ can be extended to $L^*$ with values in $\Q$ and
the \emph{discriminant quadratic form} of $L$ is the map defined by
$$q_L:A_L\lra \Q/2\Z,\ \ \ q_L(x+L)=(x,x)+2\Z.$$
\textbf{Examples:}\ \\
i) The \emph{hyperbolic plane} is the rank two even lattice
$U=\langle x,y\rangle$ with bilinear form: $(x,x)=0, (x,y)=1$.
This is an even unimodular lattice of signature $(1,1)$.\\
ii) \emph{Lattices of type $ADE$.} We denote by $A_m, D_n, E_l$,
$(m\geq 1,\ n \geq 4,\ l=6,7,8)$ the negative definite even lattices
associated with the corresponding Cartan matrixes. By the
classification theorem of root systems (see \cite{B}), any negative
definite lattice generated by elements with square $-2$ is an
orthogonal sum of lattices of type $ADE$. The discriminant groups of
these lattices are:
\begin{Prop}[Proposition 3.2.2, \cite{N2}]\label{es2}
$$A_{A_m}\cong \Z/(m+1)\Z;$$
$$A_{D_n}\cong \left\{\begin{array}{ll}
\Z/4\Z, &if\ n\ is\ odd,\\
(\Z/2\Z)^2, &if\ n\ is\ even;
\end{array}\right.$$
$$A_{E_l}\cong \Z/(9-l)\Z;$$
where $m\geq 1$, $n \geq 4$, $l=6,7,8$.
\end{Prop}
An even lattice $S$ is $2$-\emph{elementary} if $A_S\cong
(\Z/2\Z)^{\alpha_S}$ for some non negative integer $\alpha_S$.
We can associate to $S$ the following invariant:
$$\delta_S=\left\{ \begin{array}{ll}
0 & \mbox{ if }q_S(x)\in\Z \mbox{ for all }x\in A_S\\
1 & \mbox{ otherwise}
\end{array}\right.$$
The following result gives an answer to the problem of classification of $2$-elementary lattices:
\begin{Thm}[Theorem 3.6.2, \cite{N1}]\label{2el}
The isomorphism class of an even 2-elementa\\ry lattice $S$ is
determined by the invariants $(s_+,s_-,\alpha_S,\delta_S)$ if
$s_+>0$ and $s_->0$.
\end{Thm}
Examples of $2$-elementary lattices are given by unimodular lattices, ``twisted" unimodular lattices (e.g. $U(2)$, $E_8(2)$), $A_1$ and $D_{2n}$, $n\geq 2$.
We describe in detail two cases:\\
i) Let $A_1=\langle x\rangle$, then $A_{A_1}=\langle x/2\rangle\cong\Z_2$ and $q_{A_1}(x/2)=-1/2$. In particular $\alpha_{A_1}=\delta_{A_1}=1$. This lattice has a (unique up to isometries) structure of free $\Z[i]$ module given by the action of the isometry:
$$J_1=\left(\begin{array}{cc}
0 &1\\
-1 & 0
\end{array}\right).$$
ii) The lattice $D_4$ can be described as:
$$D_4=<\pm e_i\pm e_j>,\  i,j=1,\dots,4,$$
where $(e_i\cdot e_j)=-\delta_{ij}$.
In fact, a basis is given by:
$$f_1=e_1-e_2,\ f_2=e_3-e_4,\ f_3=e_2-e_3,\ f_4=e_3+e_4.$$
We have
$$A_{D_4}= \langle e_1, (e_1+e_2+e_3+e_4)/2\rangle=\langle (f_2+f_4)/2, (f_1+f_2)/2 \rangle$$
(where we write representatives of the classes).
The quadratic form, with respect to the given basis for $A_{D_4}$ is given by:
$$\left(\begin{array}{lr}
1 & 1/2\\
1/2 & 1
\end{array}\right).$$
In particular $\alpha_{D_4}=2$ and $\delta_{D_4}=0$.
The lattice $D_4$ is a free $\Z[i]$-module of rank two given by the action of the isometry:
$$J_2=\left(\begin{array}{cccc}
0  &1  &0 &0\\
-1 & 0 & 0 & 0\\
1& 0 & 1& 1\\
-1 & 1 &2 &1
\end{array}\right)$$
with respect to the basis $f_1,\dots,f_4$ (it can be proved that this structure is unique
up to isometries). More easily, the action on $e_1,e_3,e_2,e_4$ is given by the matrix $J_1\oplus J_1$.
\begin{Rem}\label{d4}
Let $r=a_1e_1+\ldots +a_4e_4$,\ $\sum a_i \equiv 0\, (mod\, 2)$ be a primitive vector in $D_4$.
If $r$ is primitive in $D_4$ (in particular not all $a_i$'s are even) we have the cases:\\
i) If $a_i$ is odd for $i=1,\dots,4$ then $r/2$ is non-trivial in $A_{D_4}$. \\
ii) If $a_1, a_3$ are odd and $a_2, a_4$ are even, then $(r +\sigma^*(r))/2$
has integer coefficients but either 1 or 3 of the
coefficients are odd, hence it is a non-trivial element in $A_{D_4}$.\\
iii) If $a_1, a_2$ are odd and $a_3, a_4$ are even, then all coefficients of $(r +
\sigma^*(r))$ are odd integers, so $(r +\sigma^*(r))/2$ is non-trivial in $A_{D_4}$.
\end{Rem}


An embedding of lattices $S\subset L$ is \emph{primitive} if $L/S$ is a free $\Z$-module. Two primitive embeddings $S\subset L$ and $S\subset L'$ are \emph{isomorphic}
if there exists an isomorphism between $L$ and $L'$ which is the identity on $S$.

Let $S\subset L$ be a primitive embedding of even lattices and $K=S^{\perp}_L$ be the orthogonal lattice to $S$ in $L$. Then $S\oplus K\subset L$ and
$$H_L=L/(S\oplus K)\subset A_S\oplus A_K$$
is isotropic. Moreover:
$$q_L=(q_S\oplus q_K)\mid{H_L}^{\perp}/H_L.$$
Let $p_S$, $p_K$ be respectively the projections of $H_L$ to $A_S$ and $A_K$. By the primitivity of the embedding, these projections are embeddings:
$$p_S:H_L\lra H_L(S)\subset A_S,\ p_K:H_L\lra H_L(K)\subset A_K.$$
The group $H_L$ is the graph of the isomorphism
$$\gamma_{S,K}^L=p_K\circ p_S^{-1}:H_L(S)\lra H_L(K)$$
in $A_S\oplus A_K$.
Since $H_L$ is isotropic we also have:
$$q_K\circ \gamma_{S,K}^L=-q_S.$$
This construction can be also reversed.
\begin{Prop}[Proposition 1.6.1, \cite{N2}]\label{emb}
A primitive embedding of an even lattice $S$ into another even lattice $L$ with discriminant form $q_L$  and orthogonal complement isomorphic to $K$ is determined by a pair $(H,\gamma)$, where:
\begin{itemize}
\item[i)] $H\subset A_S$ and $\gamma:H\lra A_K$ is a group monomorphism,
\item[ii)] $q_K\circ \gamma=-q_s\mid H$,
\item[iii)] $q_L\cong (q_S\oplus q_K\mid \Gamma^{\perp}/\Gamma)$ where $\Gamma$ is the graph of $\gamma$ in $A_S\oplus A_K$.
\end{itemize}
Two such pairs $(H,\gamma)$ and $(H',\gamma')$ determine isomorphic primitive embeddings if and only if $H=H'$ and $\gamma,\gamma'$ are conjugate via a self-isometry of $K$.
\end{Prop}
In particular, if $L$ is unimodular, the embedding is determined by an isomorphism $\gamma: A_S\lra A_K$ satisfying ii).

\end{document}